 
\baselineskip=14pt
\parskip=10pt
\def\Tilde{\char126\relax}

\font\eightrm=cmr8  
\font\eighttt=cmtt8
\magnification=\magstephalf

\parindent=0pt
\overfullrule=0in
\bf
\centerline{\bf A WZ PROOF OF RAMANUJAN'S FORMULA FOR $\pi$ }
\medskip
\it
\centerline{Shalosh B. EKHAD$^1$ and Doron ZEILBERGER\footnote{$^1$}
{\baselineskip=9pt
\eightrm  \raggedright
Department of Mathematics, Temple University,
Philadelphia, PA 19122.
{\baselineskip=9pt
\eighttt [ekhad,zeilberg]@math.temple.edu ; 
\break http://www.math.temple.edu/\Tilde 
[ekhad,zeilberg]. }
The work of the second author was
supported in part by the NSF. This paper was published
in p.107-108 of `Geometry, Analysis, and Mechanics',
ed. by J. M. Rassias, World Scientific, Singapore 1994.
}}
\medskip
\rm
\centerline{\it Dedicated to Archimedes on his $2300^{th}$ birthday}
\bigskip
 
Archimedes computed $\pi$ very accurately. Much
later, Ramanujan discovered several infinite series
for $1/ \pi$ that enables one to compute $\pi$ even more accurately.
The most impressive one is([Ra]):
($(a)_k$ denotes, as usual, $a(a+1) ... (a+k-1)$.)
 
$$
{{1} \over {\pi}} =
2 \sqrt{2} \sum_{k=0}^{\infty}
{{(1/4)_k (1/2)_k (3/4)_k}
\over {k!^3}}
(1103+26390k) (1/99)^{4k+2}.
\eqno(1)
$$
 
This formula is an example of a {\it non-terminating} 
hypergeometric series identity.
Many times, non-terminating series are either limiting cases
or "analytic continuations" of {\it terminating identities}, which
are now known to be routinely provable by computer. [WZ].
 
While we do not know of a terminating generalization of (1),
we do know how to give a WZ proof of another formula for $\pi$, also
given by Ramanujan[Ra], and included in his famous letter to
Hardy. This formula is:
$$
{{2} \over {\pi}} =
\sum_{k=0}^{\infty} (-1)^k (4k+1){{ (1/2)_k^3} \over
{k!^3}} ~~~.
\eqno(2)
$$
The terminating version, that we will prove is
$$
{{\Gamma(3/2+n)} \over {\Gamma(3/2) \Gamma ( n+1) }} =
\sum_{k=0}^{\infty} (-1)^k (4k+1) {{ (1/2)_k^2 (-n)_k
} \over
{k!^2 (3/2 + n)_k }}.
\eqno(3)
$$
 
To prove it for all {\it positive} integers $n$, we call
the summand divided by the left side $F(n,k)$, and cleverly
construct
$$
G(n,k) \,:=\,{{(2k+1)^2} \over
          {(2n+2k+3)(4k+1)}} F(n,k) \,,
$$
with the motive that $F(n+1,k)-F(n,k)=G(n,k)-G(n,k-1)$
(check!), and summing this last identity w.r.t $k$ shows that
$\sum_k F(n,k) \equiv Constant$, which is seen to be $1$, by plugging
in $n=0$. This proves (3). To deduce (2), we "plug" in
$n= -1/2$, which is legitimate in view of Carlson's theorem [Ba].
 
{\bf REFERENCES}
 
[Ba] W.N. Bailey, ``Generalized Hypergeometric Series'', (Cambridge
Univ. Press, 1935), p. 39.
 
[Ra] K.S. Rao, in  ``Srinivasa Ramanujan'', ed. K.R. Nagarajan and
T. Soundararajan, (MACMILLAN INDIA, Madras, 1988).
 
[WZ] H. S. Wilf and D. Zeilberger, {\it Amer. Math.Soc.} {\bf B3} 
(1990) 147.

\bye